**Кіяновська Наталія Михайлівна**
кандидат педагогічних наук, старший викладач кафедри вищої математики
Державний вищий навчальний заклад «Криворізький національний університет», м. Кривий Ріг, Україна
*kiianovska.nataliia@yandex.ru*

**Рашевська Наталя Василівна**
доцент, кандидат педагогічних наук, доцент кафедри вищої математики
Державний вищий навчальний заклад «Криворізький національний університет», м. Кривий Ріг, Україна
*nvr1701@gmail.com*

**Семеріков Сергій Олексійович**
професор, доктор педагогічних наук, завідувач кафедри фундаментальних і соціально-гуманітарних дисциплін
Державний вищий навчальний заклад «Криворізький національний університет», м. Кривий Ріг, Україна
*semerikov@gmail.com*

# ЕТАПИ РОЗВИТКУ ТЕОРІЇ І МЕТОДИКИ ВИКОРИСТАННЯ ІНФОРМАЦІЙНО-КОМУНІКАЦІЙНИХ ТЕХНОЛОГІЙ У НАВЧАННІ ВИЩОЇ МАТЕМАТИКИ СТУДЕНТІВ ІНЖЕНЕРНИХ СПЕЦІАЛЬНОСТЕЙ У СПОЛУЧЕНИХ ШТАТАХ АМЕРИКИ

**Анотація.** У статті досліджуються проблеми розвитку інформаційно-комунікаційних технологій (ІКТ) навчання вищої математики студентів інженерних спеціальностей у США. Охарактеризовано сутність конвергенції тенденції інформатизації системи вищої інженерної освіти США з іншими тенденціями її розвитку; визначено основні історико-педагогічні етапи розвитку теорії і методики використання ІКТ у навчанні вищої математики студентів інженерних спеціальностей у США. Вивчення історико-педагогічної джерельної бази надало можливість виокремити шість етапів, на ранніх етапах проаналізовано провідні засоби ІКТ навчання вищої математики, вказано протиріччя і визначено основні особливості використання засобів ІКТ у навчанні вищої математики студентів інженерних спеціальностей.

**Ключові слова:** досвід США; ІКТ навчання; ІКТ в освіті; ІКТ навчання вищої математики.

## 1. ВСТУП

На сучасному етапі розвитку інформаційного суспільства використання засобів інформаційно-комунікаційних технологій (ІКТ) сприяє глобалізації освіти, розвитку міжнародного ринку праці, зростанню різних видів мобільності особистості. Важливим наслідком глобалізації є підвищення мобільності студентів, абітурієнтів та випускників університетів. Зростання академічної мобільності, уведення міжнародних норм і стандартів, за допомогою яких академічні кваліфікації з різних країн можуть бути порівняні й визнані, призводить до збільшення конкуренції між ВНЗ і сприяє підвищенню якості вищої освіти.

Необхідною умовою суспільного й економічного розвитку будь-якої країни є інвестиції в освіту населення. У цьому контексті глобалізація освіти сприяє особистісному і професійному розвитку фахівців, які займаються розробкою й упровадженням нових технологій, — інженерів.

Вищі технічні навчальні заклади США мають значні педагогічні досягнення і розвинену систему підготовки фахівців інженерних напрямів на основі системного використання засобів ІКТ. У глобалізованому просторі вищої освіти проблему підвищення якості підготовки фахівців у вітчизняних ВНЗ доцільно розв'язувати через інтеграцію з кращими здобутками світової педагогічної думки і творче використання досвіду передових ВНЗ інженерного профілю.

**Постановка проблеми.** Використання ІКТ у процесі навчання вищої математики студентів інженерних спеціальностей створює умови для самореалізації студента, що сприяє підвищенню його пізнавальної активності, розвитку критичного мислення, формуванню у студентів навичок організації самостійної роботи, розвитку творчих здібностей, підвищенню відповідальності за результати своєї праці, а також вдосконаленню процесу навчання.

Тому виникає необхідність дослідження історії та сучасного стану розвитку засобів ІКТ навчання вищої математики студентів інженерних спеціальностей у технічних ВНЗ США, що займають найвищі позиції у рейтингу найкращих ВНЗ світу [38], з метою модернізації системи вищої інженерної освіти України і її спрямування на підготовку фахівців, здатних до швидкого просування науково-технічного прогресу.

**Аналіз останніх досліджень і публікацій.** У роботах В. Ю. Бикова, Ю. Г. Запорожченко, М. П. Лещенко, О. М. Спіріна, О. В. Овчарук, Н. В. Сороко, Б. І. Шуневича та інших здійснені порівняльно-педагогічні дослідження щодо зарубіжного досвіду застосування ІКТ в освіті. Теорія і методики використання ІКТ у навчанні вищої математики розроблялись у роботах К. В. Власенко, Ю. В. Горошка, В. І. Клочка, С. А. Ракова, О. В. Співаковського та інших.

Існують певні суперечності, зокрема, між сучасними вимогами до фахівця інженера і реальним рівнем їх підготовки у ВНЗ, прагненням підвищувати кваліфікацію викладачів математичних дисциплін і рівнем їх обізнаності у засобах ІКТ навчання математичних дисциплін. Залишаються недослідженими загальні тенденції розвитку засобів ІКТ навчання вищої математики студентів інженерних спеціальностей у США у контексті їх еволюції і конвергенції.

**Мета статті.** З огляду на це, **метою** статті є здійснення аналізу процесу розвитку засобів ІКТ навчання вищої математики студентів інженерних спеціальностей вищих технічних навчальних закладів у США і виділення етапів розвитку теорії і методики використання засобів ІКТ у навчанні вищої математики студентів інженерних спеціальностей у США.

## 2. МЕТОДИ ДОСЛІДЖЕННЯ

Під час дослідження використовувались такі методи: аналіз науково-педагогічної, методичної, історичної літератури, законодавчої бази США і міжнародних організацій у галузі освіти з метою встановлення стану впровадження засобів ІКТ в освітній практиці технічних ВНЗ США; аналіз, синтез та теоретичне узагальнення досвіду США з використання ІКТ у процесі навчання вищої математики; історико-педагогічний аналіз літератури для встановлення хронологічних меж основних етапів розвитку теорії і методики використання ІКТ в освітній практиці технічних ВНЗ США; моделювання для з'ясування особливостей і узагальнення досвіду використання ІКТ у навчанні інженерів у США.

## 3. РЕЗУЛЬТАТИ ДОСЛІДЖЕННЯ

На сьогоднішній день розроблено достатню кількість програмних засобів, що можуть бути використані у процесі навчання вищої математики. Деякі з них уже органічно вписалися в педагогічний процес, а для деяких розроблення методики використання знаходиться на етапі становлення.

Н. Сінклер (NathalieSinclair) [22, с. 235–253] пропонує класифікувати використання ІКТ навчального призначення не за математичним змістом, а за способом їх застосування.

Проведений історико-педагогічний аналіз літератури щодо впровадження ІКТ у навчання вищої математики і розвиток теорії і методики використання ІКТ у навчанні вищої математики студентів інженерних спеціальностей у США надав можливість виокремити такі етапи:

– *перший етап* — 1965–1973 рр. Нижня межа етапу пов'язана з появою мінікомп'ютера PDP-8, вартість якого була значно нижчою за його попередників, що призвело до закупівлі PDP-8 для навчальних цілей. Головна тенденція етапу пов'язана з появою у ВНЗ США достатньої кількості комп'ютерних засобів, що були оснащені мовами високого рівня, і специфікою апаратного забезпечення ІКТ (використання мейнфреймів з обмеженим мережним доступом);

– *другий етап* — 1973–1981 рр. Нижня межа етапу пов'язана з проведенням в Університеті Північної Кароліни IV-ої Конференції з використання комп'ютерів у навчальному процесі, на якій було відзначено потенціал UNIX для навчального процесу. Тенденція етапу пов'язана з поширенням в університетах США мережної операційної системи UNIX, використанням міні- і мікрокомп'ютерних систем;

– *третій етап* — 1981–1989 рр. Нижня межа етапу відповідає появі ОС MS DOS і першого масового використання персонального комп'ютера IBM PC. Тенденція етапу пов'язана з поширенням персональних комп'ютерів для навчання вищої математики у технічних ВНЗ США;

– *четвертий етап* — 1989–1997 рр. Нижня межа етапу пов'язана з винаходом Тімом Бернерс-Лі (Tim Berners-Lee) WorldWideWeb. Тенденція етапу пов'язана з використанням технологій Web 1.0 у навчанні вищої математики в технічних ВНЗ США;

– *п'ятий етап* — 1997–2003 рр. Нижня межа етапу пов'язана з виникненням LMS і появою нормативного документа уряду США щодо створення стандартів у сфері навчання за допомогою мережі Інтернет і мультимедіа. Тенденція етапу пов'язана з упровадженням систем управління навчанням у процес навчання вищої математики;

– *шостий етап* — з 2003 р. до теперішнього часу. Нижня межа етапу пов'язана з визнанням більшістю ВНЗ США масових відкритих дистанційних курсів й офіційним відкриттям у Масачусетському технологічному інституті сайту OpenCourseWare. Тенденція етапу пов'язана з перенесенням у Web-середовище засобів підтримки математичної діяльності і становленням і розвитком хмарних технологій навчання.

Пропонуємо почати опис з ранніх етапів (перших трьох).

**Перший етап** — 1965–1973 рр. До 1965 року використання ІКТ у навчанні вищої математики не мало системного характеру, незважаючи на вдалі, але поодинокі спроби. Показовим прикладом використання ІКТ у навчанні вищої математики в США є оголошення в липневому номері журналу «NewScientist» за 1957 рік, у якому серед кваліфікаційних вимог до посади викладача (помічника лектора — асистента) були знання мов програмування (насамперед FORTRAN) і комп'ютерної техніки.

До речі, назва цієї мови походить від FORmulaTRANslation («переклад формул» мовою, зрозумілою для комп'ютера) і відображає прикладний аспект інженерної математичної підготовки. Проте пряме перенесення програмування цією мовою у навчальний процес вищої школи зазнало утруднень, пов'язаних із:

1) станом розвитку засобів ІКТ: формування мейнфреймів — «великих» комп'ютерів високої вартості, — обслуговування яких було спрямоване на зменшення витрат від пристроїв, тому провідним режимом роботи таких комп'ютерів був пакетний (неінтерактивний режим виконання певної послідовності програм з відкладеним аналізом результатів їх роботи), у той час як для навчання провідним режимом роботи мав бути діалоговий;
2) обмеженими фінансовими можливостями закладів освіти, наслідком яких було використання насамперед застарілих мікрокомп'ютерів, що часто не мали засобів розробки мовою FORTRAN;
3) нерозробленістю психолого-педагогічних засад використання засобів ІКТ у навчанні.

На розв'язання останньої з проблем й була спрямована робота Б. Скіннера і Н. А. Кроудера (NormanAllisonCrowder) [23] з так званого «програмованого навчання». У сучасній зарубіжній психології запропонований ними підхід носить назву «інструкціоністського навчання» — гілки біхевіоризму.

У концепції програмованого навчання передбачалась така організація процесу засвоєння знань, умінь і навичок, що на кожному етапі навчального процесу чітко обумовлювались ті знання, уміння і навички, що мають бути засвоєні, і контролювався процес їх засвоєння [1].

Головна ідея цієї концепції – управління учінням, навчально-пізнавальними діями студентів за допомогою навчальної програми.

Розв'язання другої проблеми було виконано у 1964 році Дж. Кемені (John Kemeny) та Т. Курцем (Thomas Kurtz) [1819]. Створена ними мова BASIC була чи не найпершою спробою з реалізації інструкціоністського навчання на рівні мови програмування.

Не слід вважати, що лише комп'ютери у вузькому сенсі були об'єктом інформатизації в системі освіти США — так, у 1964 році відбулась публічна дискусія у коледжі Хоупа на тему надання високошвидкісних комп'ютерів й електронних калькуляторів для потреб навчання [27].

Електронні калькулятори (за сутністю спеціалізовані мікрокомп'ютери) у навчанні вищої математики у той час були провідними, проте далеко не єдиними засобами — так, у [24] наводиться перелік засобів ІКТ для навчання математичної логіки: комп'ютерний термінал з можливістю візуального й аудіального подання навчальних матеріалів, клавіатура для введення письмових відповідей, мікрофон для аудіовідповідей і світлове перо для вибору об'єктів.

Отже, на початок 1965 року у системі освіти США була достатня кількість комп'ютерних засобів різного рівня, оснащених мовами високого рівня, що надає можливість вважати цей рік умовною нижньою межею першого етапу розвитку теорії і методики використання ІКТ у навчанні вищої математики студентів інженерних спеціальностей у США.

У 1965 році фірмою DEC (Digital Equipment Corporation) був випущений перший комерційно успішний мінікомп'ютер — PDP-8.

Як згадується у [11; 13], комп'ютер PDP-8 використовували на заняттях з математики на математичних факультетах. У статті [13] показано, як два вчителі математики запропонували створення окремого комп'ютерного відділу, частково, для навчання дітей, які не мають математичних схильностей, і частково тому, що вони не встигали надавати консультаційну допомогу учням за комп'ютером, якщо вони цього потребували, одночасно виконуючи іншу роботу.

Сімейство мінікомп'ютерів PDP часто зустрічається у наукових розвідках з історії ІКТ саме через їх поширеність. Оптимальне на той час поєднання засобів ІКТ за поміркованих цін сприяло поширенню цього сімейства у наукових і освітніх установах. Одним із дослідницьких проектів зі створення діалогової математичної системи стала Reduce (www.reduce-algebra.com) — система комп'ютерної алгебри загального призначення.

Дослідження Ж. Піаже (Jean William Fritz Piaget) із психології раннього дитинства надали можливість його учню С. Пейперту (Seymour Papert) у 1967 році у Массачусетському технологічному інституті, не відходячи повністю від інструкціонізму, запропонувати новий засіб навчання — мову LOGO [1819], що базувався на конструктивістському підході до навчальної діяльності. У середовищі LOGO її користувач — програміст — виступав у ролі «вчителя» для головного об'єкта мікросвіту LOGO — черепахи, «навчаючи» її через програмування виконувати певні дії.

Як зазначено в [19], мова програмування LOGO призначена для заохочування до строгого мислення в математиці. Дизайн середовища LOGO справив значний вплив на подальший розвиток засобів навчання і навчальні концепції. Так, один із співробітників С. Пейперта — А. Кей (Alan Curtis Kay) запропонував у 1968 році Dynabook [6] — «персональний комп'ютер для дітей будь-якого віку», оснащений середовищем мови Smalltalk (www.smalltalk.org/main).

Це була перша мова високого рівня, що підтримувала експериментування з широким набором математичних об'єктів — від чисел з різних множин до геометричних об'єктів.

Отже, з виходом LOGO і Smalltak-72 програмоване навчання перестало бути домінуючою концепцією, що зумовило вибір верхньої межі першого етапу.

Знову звернувшись за приклад до історії використання ІКТ у коледжі Хоупа, можна прослідкувати їх еволюцію протягом першого етапу:

1967 рік — створення обчислювальної лабораторії й обладнання її калькуляторами IME-86 [14];

1968 рік — застосування IBM 1130 для генерування псевдовипадкових чисел [2];

1969 рік — проект Національного наукового фонду США «Використання комп'ютерів у навчанні статистики» [26];

1970 рік — дворічний інтегрований курс «Прикладна статистика та програмування» мовою FORTFAN [27];

1972 рік — видане керівництво з виконання лабораторних робіт з теорії ймовірностей і математичної статистики [30].

Традиційно не залишилися осторонь й провідні ВНЗ: так, у МТІ у 1968 році у рамках проекту MAC (Mathematics and Computation) було створено систему комп'ютерної математики Macsyma (maxima.sourceforge.net), а у 1971 році IBM під керівництвом Р. Дженкса (Richard Jenks) — систему Axiom (www.axiom-developer.org), що в той час мала назву Scratchpad.

Винайдення у 1971 році флоппі-диска сприяло персоналізації використання мейнфреймів і мінікомп'ютерів: 81,6 Кб даних, що уміщувались на 8-дюймовому диску, було цілком достатньо для зберігання текстових документів (статей, програм тощо), найбільш поширених у академічному середовищі [32].

Для підтримки навчання математики у цей період було розроблено низку спеціалізованих пристроїв:

– KENBAK-1 (розробник — Дж. Бланкенбейнер (JohnBlankenbaker), 1972 рік) — перший комп'ютер для навчання, розрахований на непрофесійних користувачів. Як зазначається у [3], цей комп'ютер поєднував гнучкість і доступність;

– HP-35 (розробник — Hewlett Packard, 1972 рік) —один з перших калькуляторів HP, що містив мікропроцесор. Це сприяло його компактності і зручності використання у навчанні математики. Дж. К. Хорн (Joseph K. Horn), який використовував калькулятори HP як учитель математики [15], є автором багатьох статей з їх застосування. Як спеціалізований мікрокомп'ютер, HP-71 (подальший розвиток HP-35) надавав можливість програмування мовою BASIC, а його наступник — мовою СКМ Derive.

Як зазначає Дж. Г. Харві (John G. Harvey) [12], у навчанні вищої математики калькулятори доцільно застосовувати, зокрема, для тестування.

HP-35 був розроблений «для інженерів та студентів інженерних спеціальностей» [12, с. 140]. Наукові і графічні калькулятори HP (як програмовані, так й непрограмовані) підтримують й нетрадиційні для калькуляторів дії над матрицями і статистичні функції.

Головною проблемою, на думку Дж. Г. Харві [12, с. 145], була проблема нерозробленості методики ефективного використання калькуляторів у навчанні математики, з одного боку, і нерозробленості методики навчання математики, орієнтованої на використання калькуляторів.

Протекціонізм комп'ютерних фірм до окремих ВНЗ привів до суттєвого зростання забезпеченості навчальних закладів засобами ІКТ: так, якщо у 1965 році менше 5 % усіх навчальних закладів були забезпечені комп'ютерами для навчальних потреб [31, с. 51], то у 1972 році передавання даних через мережу набуло актуальності за рахунок суттєвого зростання забезпеченості комп'ютерною технікою і комунікаційними засобами [36, с. 11].

Проведений аналіз надає можливість визначити основні особливості використання засобів ІКТ у навчанні вищої математики студентів інженерних спеціальностей на першому етапі їх розвитку:
1) діалоговий режим роботи з навчальними програмами;
2) поява перших систем підтримки математичної діяльності без програмування мовами загального призначення;
3) розмаїття апаратного і програмного забезпечення;
4) домінування біхевіоризму в обґрунтуванні використання ІКТ і розробці навчальних програм;
5) уведення програмування в курси вищої математики.

Указані особливості породили низку протиріч:
1) між недостатньою адекватністю програмованого навчання для опису навчальної діяльності й особистим розвитком студента у процесі навчання;
2) між необхідністю посилення прикладного аспекту навчання вищої математики майбутніх інженерів і недостатністю навчального часу на одночасне опанування математики і програмування;
3) між розмаїттям засобів ІКТ і необхідністю уніфікації засобів навчання.

Часткове розв'язання вказаних протиріч було досягнуто ще на першому етапі — з'явились нові, більш адекватні підходи до моделювання процесу навчання, перші системи комп'ютерної математики й у 1969 році — операційна система UNIX [17], спрямована на об'єднання різних засобів ІКТ у єдиному мережному середовищі.

**Другий етап** — 1973–1981 рр. — пов'язаний із поширенням в університетах США мережної операційної системи UNIX, використанням міні- і мікрокомп'ютерних систем.

Із спогадів викладача вищої математики з Каліфорнійського університету в Берклі [13]: «було вирішено створити середовище в середній школі, що було б найбільш близьким до лабораторій МТІ і Стенфордського університету, де навчався я.Було встановлено на PDP-11/70 версію UNIX 7; ми були альфа-тестувальниками BSD 2.9, версії BerkeleyUnix для PDP-11. Установка, тестування і налагодження цієї нової системи було проведено виключно за рахунок студентів».

Цей фрагмент містить згадку про PDP-11 — подальший розвиток застосовуваної на першому етапі PDP-8 і показує зацікавленість фірми-виробника (DEC) у наданні власних засобів ІКТ університетам з ОС UNIX. Версія UNIX — BSD — відображає внесок університету в її розробку (Berkeley Software Distribution) [33].

Головна особливість UNIX — мобільність (насамперед, програмна) — створила умови для поширення програмного забезпечення під управлінням цієї системи на різні засоби ІКТ — від мейнфреймів до міні- (а надалі і мікро-) комп'ютерів.

Багато відомих на теперішній час систем комп'ютерної математики були створені саме у цій ОС: MATLAB (наприкінці 1970-х), Maple (1979) та інші були розроблені як мови програмування для навчання студентів вищої математики, використовуючи різні математичні бібліотеки, не вивчаючи мову FORTRAN.

У зв'язку з тим, що ОС UNIX і її програмне забезпечення були мобільними, з'явилась можливість для об'єднання не лише обчислювальних ресурсів різних комп'ютерів, а й користувачів, їх програм і даних у мережному середовищі: «отримана система UNIX надавала користувачам можливості віддаленого доступу до терміналу і спільно використовувати файлові системи; вихідний код поставлявся із системою, користувачі могли обмінюватися даними і програмами безпосередньо і неофіційно; оскільки UNIX працювала на відносно недорогій міні-ЕОМ, малі групи дослідників могли вільно експериментувати з нею» [35].

Потенціал UNIX для навчального процесу був відзначений на IV Конференції з використання комп'ютерів у навчальному процесі, що відбулася в 1973 році в Університеті Північної Кароліни [27].

Початок другого етапу відзначався узагальненням досвіду використання засобів ІКТ навчання вищої математики; зокрема, на восьмій конференції з використання комп'ютерів у навчанні студентів [25] (1977), конференції НАТО з комп'ютерно-орієнтованого навчання [27] (1976) та інших.

У середині 1970-х рр. з'явився Журнал мічиганської асоціації користувачів комп'ютера у навчанні, у якому в 1978 році було опубліковано фундаментальну роботу з навчання теорії ймовірностей і математичної статистики із застосуванням комп'ютерів [28].

Розпочата після винаходу дискет персоналізація засобів ІКТ була продовжена розробками 1977 року — комп'ютерами Apple II і TRS-80.

Застосування конструктивістського підходу до навчання вищої математики сприяли нові засоби Apple II для пересічних користувачів: графічний дисплей і маніпулятор «миша». Це надавало можливість створювати навчальні ігри, системи комп'ютерної графіки і динамічної геометрії. А вбудована мова Apple II — BASIC — була доповнена графічними командами. Конструкція виявилась настільки вдалою, що стала родоначальником не лише комп'ютерів Apple Macintosh, а й вітчизняного комп'ютера для системи освіти «Агат» [4].

Як зазначено в [5], у рамках конструктивізму, методичною метою розробленого проекту навчання математичних дисциплін на комп'ютерах Apple II було створення навчального середовища, що сприяло набуттю студентами математичних понять через повторювані цикли розроблених завдань, вивчаючи завдання і розмірковуючи про способи їх розв'язання.

У [20] згадується про використання комп'ютерів TRS-80 як засобів навчання в аудиторії. Викладачі використовували програму на Level II BASIC для комп'ютера TRS-80, що імітувала машину Тюрінга і демонструвала природу пристрою. Програма запускалася в динамічному режимі і була призначена для використання як навчального посібника з інформатики або математики для студентів, які вивчали теорію числення або теорію автоматів.

Поширення персоналізованих комп'ютерів з графічним інтерфейсом, орієнтованих на ігрову діяльність (Atari 800 та інші) також сприяло розвитку ігрових засобів навчання математики [7].

Опубліковані у 1978 році «Рекомендації з виконання лабораторних робіт з теорії ймовірностей та математичної статистики» [29] зіграло значну роль у розвитку теорії і методики навчання математичних дисциплін із застосуванням комп'ютерів. Через 15 років після перших спроб уведення комп'ютерів у навчання вищої математики майбутніх інженерів відчувається суттєве зміщення: з навчання математики разом з програмуванням до обґрунтованого використання програмних засобів у процесі навчання.

Індустрія програмних засобів навчального призначення на другому етапі була орієнтована переважно на Apple II [10].

У 1975 році інженерами Xerox PARC (тієї самої лабораторії, де були винайдені миша, графічний інтерфейс і Smalltalk) був отриманий патент США на технологію Ethernet [19], а вже у 1979 році з'явилась перша відкрита інформаційна онлайн-служба Compuserve. І хоча ядром її були UNIX-системи, об'єднувала вона й персональні комп'ютери.

Отже, можна виокремити такі характерні риси другого етапу розвитку теорії і методики використання ІКТ у навчанні вищої математики студентів інженерних спеціальностей у США:
1) перехід від використання мов програмування у навчанні до використання математичних бібліотек, систем комп'ютерної математики і мов високого рівня;
2) застосування комп'ютерної графіки у навчальних програмах;
3) поява нових класів навчальних програм — навчальних ігор, систем динамічної геометрії й електронних таблиць;
4) поява і поширення комп'ютерних мереж, що сприяли активному спілкуванню між викладачами і студентами;
5) розвиток засобів ІКТ навчання вищої математики — графічних і символьних калькуляторів.

У розвитку засобів ІКТ на цьому етапі можна виділити низку протиріч:
1) між потенціалом використання мультимедійних засобів комп'ютерних систем і не розробленістю психолого-педагогічних основ їх використання;
2) між потребою студентів і викладачів у персональних засобах ІКТ та недостатністю пропозицій виробників;
3) між потребою виробників персональних комп'ютерів у адаптованому для них варіанті UNIX і недостатнім апаратним забезпеченням її функціонування.

Вказані протиріччя визначили початок третього етапу, на якому вони були вирішені.

**Третій етап** — 1981–1989 рр. — пов'язаний із поширенням персональних комп'ютерів.

Нижня межа етапу (1981 рік) відповідає появі ОС MS DOS і комп'ютера IBM PC, що відіграли надзвичайно значну роль, сформувавши сучасний ринок програмних засобів навчання вищої математики.

Персоналізація (користувача, комп'ютера, програм, даних) була провідною концепцією цього етапу. Навчання з використанням комп'ютера є одним з найбільш поширених підходів, призначених підвищити навчальні досягнення студентів вищої школи.

У ці роки значна кількість робіт була присвячена психолого-педагогічному обґрунтуванню навчання з використанням персональних комп'ютерів і розробленню навчальних курсів з дисциплін. Аналізуючи публікації того часу, можна зазначити, що найбільш поширеними мовами програмування були BASIC, Pascal, FORTRAN, Algol. Так, у Массачусетському технологічному інституті на заняттях з вищої математики студенти займались написанням комп'ютерних програм для проведення досліджень з різних тем математики. Тенденції і досягнення в галузі численних методів, обчислення і новаторські дослідження в галузі прикладної математики сприяли впровадженню комп'ютерно орієнтованих технологій навчання прикладної математики [34].

Дж. Енгельбрехт (Johann Engelbrecht) [8] виявляє низку переваг застосування комп'ютерів у процесі викладанні вищої математики: відпрацювання практичних умінь і навичок; розмаїття подання навчального матеріалу; використання у процесі моделювання і програмування.

Викладач кафедри вищої математики університету Каліфорнії Дж. Дж. Ваврік (John J. Wavrik) [37] — автор курсу «Комп'ютерна алгебра» для студентів у символьних обчисленнях (середина 80-х років ХХ століття) один із перших стверджував, що використання комп'ютерів є необхідним у процесі навчання вищої математики, оскільки їх використання сприяє підвищенню рівня навчальних досягнень студентів [16].

Стаття Дж. Дж. Вавріка [37] «Комп'ютери і кратні корені полінома» «…не описує різні повороти долі, що ведуть працівника в чистій галузі математики такій, як обчислювальна геометрія, пов'язаної з комп'ютерами. Цілком ймовірно, що персональні комп'ютери стають більш поширеними серед алгебраїстів, які роблять акцент на використання їх для науково-дослідної роботи… Природно сподіватися, що комп'ютери можуть бути використані для полегшення складності обчислень. Як виявилося, цей процес не такий уже й і простий, як спочатку може здатися… Ми хотіли б мати машину, що може допомогти, наприклад, в обчисленнях інваріантів для специфічних станів об'єктів дослідження…» [37]. У статті автор, розглядаючи питання знаходження кратних коренів полінома, зауважує, що «проблеми ефективного обчислення найбільшого спільного дільника двох многочленів з цілими коефіцієнтами отримала велику увагу… Комп'ютерні системи, призначені для алгебраїчних обчислень, надали можливість «з нескінченною точністю» робити арифметичні обчислення та використовувати для алгоритмів систем цього типу. Машини, що виконують обчислення із заданою точністю часто допускають точність лише між 6 і 16 цифрами наближень... Ця стаття містить комп'ютерну програму, що служить прикладом реалізації на конкретній машині. Програма, написана на BASIC, — найбільш поширеній мові для мікрокомп'ютерів» [37]. Автор вказує, що «метою даної статті є представлення читачеві деяких цікавих аспектів алгебраїчних обчислень. Програма, включена в статтю, не подається читачам як частина готового програмного забезпечення, а скоріше, щоб дати їм ідеї для розробки своїх власних програм» [37].

У 1984 році Р. Д. Пеа (Roy D. Pea) [21] використовував програму AlgebraLand, застосування якої надало студентам можливість автоматизувати алгебричні обчислення і зосередити увагу на розв'язанні завдань більш високого рівня складності. AlgebraLand була розроблена для того, щоб з її допомогою студенти «вивчали питання приросту швидкості», закріплюючи навички під час розв'язання задач.

Якщо звернутися за приклад до історії використання ІКТ у коледжі Хоупа, можна прослідкувати їх еволюцію протягом третього етапу [26]:

1982 рік — перша Міжнародна конференція з викладання статистики (ICOTS I) «Використання мікрокомп'ютерів для розуміння понять теорії ймовірностей та математичної статистики»;

1983 рік — щорічні збори Математичної асоціації Америки в Денвері «Використання мікрокомп'ютерів надає можливість ілюструвати основні поняття теорії ймовірностей та математичної статистики»;

1984 рік — конференція з використання мікрокомп'ютерів у статистиці Американської асоціації статистиків Університету штату Делавер «Використання мікрокомп'ютерів у галузі викладання теорії ймовірностей та математичної статистики»;

1985–1996 роки — створення обчислювальної лабораторії «Комп'ютерна лабораторія вступу до статистики» та обладнання її ПК IBM з BASIC;

1987 рік — щорічні збори Американської статистичної асоціації у Сан-Франциско, Каліфорнія «Комп'ютерне моделювання для підтвердження теоретичних понять».

Значного розвитку на цьому етапі набувають системи комп'ютерної математики, а саме: Cayley (1982), MathCAD (1985), Fermat (1985), GAP (1986), Mathematica (1986), Derive (1988), MuPad (1989).

У звіті 1989 р. Національної академії наук США про майбутнє математичної освіти [9] (так, як воно бачилось наприкінці третього етапу) вказується, що символьні комп'ютерні системи вимагають фундаментального переосмислення методики навчання математики з використанням комп'ютера: якщо до цього використовувались обчислювальні і графічні програми, що мали похибки округлення та візуального подання, то символьні обчислення є точними і відповідають діям, що виконує людина. «Пріоритети математичної освіти повинні змінюватися з метою відображення шляхів використання комп'ютерів у математиці» [9, с. 78]: провідними засобами навчання вищої математики укладачі звіту вважали електронні таблиці, пакети числового аналізу, символьні комп'ютерні системи, засоби графічного подання відомостей, перспективними — електронні підручники, віддалені класні кімнати, інтегровані навчальні середовища.

«Підручники, програмне забезпечення, комп'ютерні мережі у найближчі роки скомбінуються у новий гібридний освітньо-інформаційний ресурс» [9, с. 82] — саме цей прогноз ознаменував кінець третього і початок четвертого етапу розвитку теорії і методики використання ІКТ у навчанні вищої математики студентів інженерних спеціальностей у США.

Отже, до характерних рис третього етапу розвитку теорії і методики використання ІКТ у навчанні вищої математики студентів інженерних спеціальностей у США відносяться:

1) широке використання математичних бібліотек, систем комп'ютерної математики і проблемно орієнтованих мов;
2) широке впровадження персональних і персоналізованих засобів ІКТ у навчання математичних дисциплін;
3) використання ІКТ загального призначення (текстові редактори, електронні таблиці, бази даних тощо) для підтримки навчання математичних дисциплін.

У розвитку засобів ІКТ на цьому етапі можна виділити низку протиріч:

1) між потенціалом об'єднання використання глобальних комп'ютерних мереж і недостатньою розробленістю персональних засобів доступу до них;
2) між доцільністю перенесення гіпертекстових і гіпермедіальних систем навчального призначення у мережне середовище і недостатньою розробленістю засобів доступу до них;

3) між потребою студентів і викладачів у перенесенні навчальних матеріалів у мережу і не розробленістю психолого-педагогічних основ її використання.

Вказані протиріччя визначили початок четвертого етапу, на якому вони були вирішені.

## 4. ВИСНОВКИ ТА ПЕРСПЕКТИВИ ПОДАЛЬШИХ ДОСЛІДЖЕНЬ

Проаналізувавши джерела з проблеми дослідження, можна стверджувати, що на кожному етапі розвитку теорії і методики використання ІКТ у навчанні вищої математики студентів інженерних спеціальностей у США поступово розвивається і педагогічна наука, враховуючи науково-технічні досягнення часу. Поява нових ІКТ призводить до змін у теорії і методиці викладання вищої математики; обумовлює виникнення нових цілей, засобів, форм, методів організації процесу навчання і доповнення його змісту. Ці зміни мають позитивний вплив на педагогічний процес, значно розширюючи можливості студентів.

Перспективи подальших досліджень: провести дослідження і зробити ґрунтовний аналіз четвертого, п'ятого та шостого етапів розвитку теорії та методики використання ІКТ у навчанні вищої математики студентів інженерних спеціальностей у Сполучених Штатах Америки.

## СПИСОК ВИКОРИСТАНИХ ДЖЕРЕЛ

1. Ягупов В. В. Педагогіка : навчальний посібник / В. В. Ягупов. — К. : Либідь, 2002. — 560 с.
2. All about the IBM 1130 Computing System [online] // IBM1130.org, 2011. — Available from: http://ibm1130.org/.
3. Blankenbaker J. The first personal computer: KENBAK-1 computer [online] / John Blankenbaker // KENBAK-1 Computer, 2010. — Available from : http://www.kenbak-1.net/index.htm.
4. Bores L. D. AGAT: A Soviet Apple II Computer / Leo D. Bores //BYTE : The Small System Journal. — 1984. — Vol. 9, No. 12, November. — P. 134–135.
5. Confrey J. High school mathematics development of teacher knowledge and implementation of a problem-based mathematics curriculum using multirepresentational software :Apple classrooms of tomorrow research : Report number 11 / JereConfrey, Susan C. Piliero, Jan M. Rizzuti // Cupertino : Apple Computer, 1990. — 9 p.
6. Dalakov G. The Dynabook of Alan Kay [online] / GeorgiDalakov // History of computers, 2013. — Available from : http://history-computer.com/ModernComputer/Personal/Dynabook.html.
7. Enchin H. Home-computer programs help kids with math / Harvey Enchin // The gazette business. — Montreal. — 1983. —Aprsl 12. — P. F-2.
8. Engelbrecht J. Teaching Undergraduate Mathematics on the Internet / Johann Engelbrecht, Ansie Harding // Journal of Online Mathematics and its Applications, 2005. — №(58)2. — P. 235–276.
9. Everybody Counts: A Report to the Nation on the Future of Mathematics Education / Mathematical Sciences Education Board, Board on Mathematical Sciences, Committee on the Mathematical Sciences in the Year 2000, National Research Council. — Washington : National Academy Press, 1989. — 130 p.
10. Geometric Supposer: Triangles/64K Apple Ii, Iie, Iic/Disk, Backup, Teacher's Guide/Quick Reference Cards/Book No 1340-Fs [online]. — Cincinnati : Sunburst Communications, 1986. — Available from : http://www.amazon.com/Geometric-Supposer-Triangles-Teachers-Reference/dp/9996313344/.
11. Greenwald S. B. Computer-assisted explorations in mathematics: pedagogical adaptations across the Atlantic : a report on the development of CATAM at Cambridge University and the Mathematics Project Laboratory at MIT / Suzanne B. Greenwald, Haynes R. Miller // University Collaboration for Innovation: Lessons from the Cambridge-MIT Institute, Sense Publishers. — Massachusetts Institute of Technology, Department of Mathematics, 2007. — P. 121–131.
12. Harvey J. G. Mathematics Testing with Calculators: Ransoming the Hostages / John G. Harvey // Mathematics assessment and evaluation: Imperatives for mathematics educators / Edited by Tomas A. Romberg. — New York : State University of New York, 1992. — P. 139–168.

13. Harvey B. A case study: the Lincoln-sudbury regional high school [online] / Brian Harvey // Computer Science Division : EECS at UC Berkeley, 2012. — Available from : http://www.cs.berkeley.edu/~bh/lsrhs.html.
14. Haven G. IME-86: the only electronic desk calculator thet can grow into a compunrt [online] / Grand Haven // Hope College Department of Mathematics, 1967. — Available from : http://www.math.hope.edu/tanis/History/ime-calculator-purchase.pdf.
15. Horn J. K. The HP-71B «Math machine» 25 years old [online] / Joseph K. Horn, Richard J. Nelson, Dale Thorn // Hewlett-Packard Development Company, [2007]. — Available from :http://h71028.www7.hp.com/enterprise/downloads/The%20HP-71B%20Math%20Machine%20V2e.pdf.
16. John J Wavrik [online] / John J. Wavrik // University of California. — Available from : http://math.ucsd.edu/~jwavrik/
17. Kernighan B. W. The UNIX programming environment / Brian W. Kernighan, Rob Pike. — Prentice-Hall, 1984. — 357 p.
18. Molnar A. Computers in Education: A Brief History [online] / Andrew Molnar // THE Journal: Technological Horizons in Education. — 06.01.1997. — Available from :http://thejournal.com/Articles/1997/06/01/Computers-in-Education-A-Brief-History.aspx?Page=3.
19. Multipoint data communication system with collision detection : US 4063220 A [online] / David R. Boggs, Butler W. Lampson, Robert M. Metcalfe, Charles P. Thacker. — March 31, 1975. — Available from : http://www.google.com/patents/US4063220.
20. Navarro A. B. Turing Machine Simulator. / Navarro Aaron B. // Computers in mathematics and science teaching, 1981. — v1. — № 2. — P. 25–26.
21. Pea R. D. Cognitive technologies for mathematics education / Roy D. Pea //Cognitive science and mathematics education, 1987. — C. 89–122.
22. Sinclair N. Understanding and projecting ICT trends in mathematics education / NathalieSinclair, Nicholas Jackiw// Teaching secondary mathematics with ICT / Edited by Sue Johnston-Wilder and David Pimm. — New York : Open University Press, 2005. — P. 235–253.
23. Skinner B. F. The technology of teaching / B. F. Skinner. — New York : Appleton-Century-Crofts, 1968. — 271 p.
24. Suppes P. Computer-based mathematics instruction : The First Year of the Project : (1 September 1963 to 31 August 1964) [online] / Patrick Suppes. // Bulletin of the International Study Group for Mathematics Learnin, 1965 p. — №3. — P.7–22. — Mode of access : http://suppes-corpus.stanford.edu/article.html?id=54-1.
25. Tanis A. A computer-based laboratory for mathematical statistics and probability [online] / A. Tanis // Computer lab for statistics and probability : Computer Hope. — Available from :http://www.math.hope.edu/tanis/History/ccuc-8%201977.pdf.
26. Tanis E An Experimental Approach to the Central Limit Theorem [online] / E. Tanis, Deanna Gross // Annual Spring Meeting of the Michigan Section of the Mathematical Association of America : Hope College Department of Mathematics, 1969. — Mode of access : http://www.math.hope.edu/tanis/History/clt-deanna-1969-rotated.pdf.
27. Tanis E. Hope College : Department of Mathematics [online] / Elliot Tanis // Academic departments math. — Available from : http://www.math.hope.edu/tanis/oldpapersgiven.html.
28. Tanis E. A. Concepts in probability and statistics illustrated with the computer / Elliot A. Tanis // Michigan association of computer users for learning. —1978. —Vol II. — № 1. — P. 6-–17.
29. Tanis E. A. Laboratory manual for probability and statistical inference / Elliot A. Tanis —Iowa : CONDUIT, 1977. — 88 p.
30. Tanis E. A. Theory of Probability and statistics illustrated by the computer [online] / Elliot A. Tanis // GLCA Computing symposium at Wabash University : Hope College Department of Mathematics. – March 7–8. — 1972. — Available from : http://www.math.hope.edu/tanis/History/wabash%201972-rotated.pdf.
31. The mathematical sciences : Undergraduate education : a report // National research council (U.S.). Panel on undergraduate education in mathematics. — National Academies, 1968. — 113 p.
32. This month in ed tech history archives [online] // The Software & Information Industry Association, 2013. . — Available from :http://www.siia.net/index.php?option=com_content&view=article&id=314:this-month-archive&catid=159:education-articles&Itemid=328.
33. UNIX BSD information [online] // Computer Hope. — Available from : http://www.computerhope.com/unix/bsd.htm.
34. Usluel Y. K. A Structural Equation Model for ICT Usage in Higher Education / YaseminKocakUsluel, PetekAskar, Turgay Bas // Educational Technology & Society, 2008. — № 11 (2). — P. 262–273.
35. Vleck T. V. Dennis M. Ritchie – A. M. Turing Award Winner [online] / Tom Van Vleck // Association for Computing Machinery, 2012. — Available from : http://amturing.acm.org/award_winners/ritchie_1506389.cfm.

36. Watson P. G. Using the computer in education: A briefing for school decision makers / Paul G. Watson. — New Jersey : Educational Technology, 1972. — 128 p.
37. Wavrik J. J. Computers and the Multiplicity of Polynomial Roots / John J. Wavrik // The American Mathematical Monthly : An official journal of the Mathematical Association of America. — 1982. — Vol. 89, No. 1. — P. 34–36, 45–56.
38. World University Rankings 2010-2011 — Times Higher Education [online] // TSL Education Ltd. — London, 2012. — Available from :http://www.timeshighereducation.co.uk/world-university-rankings/2010-11/world-ranking

# ЭТАПЫ РАЗВИТИЯ ТЕОРИИ И МЕТОДИКИ ИСПОЛЬЗОВАНИЯ ИНФОРМАЦИОННО-КОММУНИКАЦИОННЫХ ТЕХНОЛОГИЙ В ОБУЧЕНИИ ВЫСШЕЙ МАТЕМАТИКИ СТУДЕНТОВ ИНЖЕНЕРНЫХ СПЕЦИАЛЬНОСТЕЙ В СОЕДИНЕННЫХ ШТАТАХ АМЕРИКИ

**Кияновская Наталия Михайловна**
кандидат педагогических наук, старший преподаватель кафедры высшей математики
Государственное высшее учебное заведение «Криворожский национальный университет», г. Кривой Рог, Украина
*kiianovska.nataliia@yandex.ru*

**Рашевская Наталья Васильевна**
доцент, кандидат педагогических наук, доцент кафедры высшей математики
Государственное высшее учебное заведение «Криворожский национальный университет», г. Кривой Рог, Украина
*nvr1701@gmail.com*

**Семериков Сергей Алексеевич**
профессор, доктор педагогических наук, заведующий кафедры фундаментальных и социально-гуманитарных дисциплин
Государственное высшее учебное заведение «Криворожский национальный университет», г. Кривой Рог, Украина
*semerikov@gmail.com*

**Аннотация.** В статье исследуются проблемы развития информационно-коммуникационных технологий (ИКТ) обучения высшей математике студентов инженерных специальностей в США. Охарактеризованы сущность конвергенции тенденции информатизации системы высшего инженерного образования США с ее другими тенденциями развития; определены основные историко-педагогические этапы развития теории и методики использования ИКТ в обучении высшей математике студентов инженерных специальностей в США. Изучение историко-педагогических источников позволило выделить шесть этапов, на каждом из которых проанализированы ведущие средства ИКТ обучения высшей математике, указано противоречия и определены основные особенности использования средств ИКТ в обучении высшей математике студентов инженерных специальностей.

**Ключевые слова:** опыт США; ИКТ обучения; ИКТ в образовании; ИКТ обучения высшей математике.

# DEVELOPMENT OF THEORY AND METHODS OF USE OF INFORMATION AND COMMUNICATION TECHNOLOGIES IN TEACHING MATHEMATICS OF ENGINEERING SPECIALITIES STUDENTS IN THE UNITED STATES

**Nataliia M. Kiianovska**
Ph.D. (in Pedagogics), AssociateProfessor of Mathematics
State Higher Educational Institution «Kryvyi Rih National University», Kryvyi Rih, Ukraine
*kiianovska.nataliia@yandex.ru*

**Natalia V. Rashevska**
Ph.D.(in Pedagogics), Associate Professor of Mathematics
State Higher Educational Institution «Kryvyi Rih National University», Kryvyi Rih, Ukraine
*nvr1701@gmail.com*

**Serghii O.Semerikov**
Professor, Doctor of Education, Head of fundamental and social and human sciences
State Higher Educational Institution «Kryvyi Rih National University», Kryvyi Rih, Ukraine
*semerikov@gmail.com*

**Abstract.** The article deals with the problems of information and communication technologies (ICT) development in teaching mathematics of engineering specialities students in the United States. In the article the nature of trends of convergence of information system in higher technical education and other tendencies in the USA are characterized. The main historical stages of development of the theory and methods of ICT use in teaching mathematics of engineering specialities students in the United States are defined. The study of historical sources has been allowed to emphasize six stages, at each stage it is analyzed the use of ICT for teaching mathematics, it is shown the contradictions and the main features of the use of ICT in teaching mathematics of engineering specialities students.

**Keywords:** american experience; ICT training; ICT in Education; ICT teaching mathematics.

# REFERENCES (TRANSLATED AND TRANSLITERATED)

1. Jaghupov V. V. Pedagogy: manual / V. V. Jaghupov. — K. : Lybidj, 2002. — 560 s. (inUkrainian)
2. All about the IBM 1130 Computing System [online] // IBM1130.org, 2011. — Available from : http://ibm1130.org/ (in English)
3. Blankenbaker J. The first personal computer: KENBAK-1 computer [online] / John Blankenbaker // KENBAK-1 Computer, 2010. — Available from : http://www.kenbak-1.net/index.htm (inEnglish)
4. Bores L. D. AGAT: A Soviet Apple II Computer / Leo D. Bores //BYTE : The Small System Journal. — 1984. — Vol. 9, No. 12, November. — P. 134–135. (in English)
5. Confrey J. High school mathematics development of teacher knowledge and implementation of a problem-based mathematics curriculum using multirepresentational software : Apple classrooms of tomorrow research : Report number 11 / JereConfrey, Susan C. Piliero, Jan M. Rizzuti // Cupertino : Apple Computer, 1990. — 9 p. (in English)
6. Dalakov G. The Dynabook of Alan Kay [online] / GeorgiDalakov // History of computers, 2013. — Available from : http://history-computer.com/ModernComputer/Personal/Dynabook.html (in English)
7. Enchin H. Home-computer programs help kids with math / Harvey Enchin // The gazette business. — Montreal. – 1983. —Aprsl 12. — P. F–2. (in English)
8. Engelbrecht J. Teaching Undergraduate Mathematics on the Internet / Johann Engelbrecht, Ansie Harding // Journal of Online Mathematics and its Applications, 2005. — №(58)2. — P. 235–276. (in English)
9. Everybody Counts: A Report to the Nation on the Future of Mathematics Education / Mathematical Sciences Education Board, Board on Mathematical Sciences, Committee on the Mathematical Sciences in the Year 2000, National Research Council. — Washington : National Academy Press, 1989. — 130 p. (in English)
10. Geometric Supposer: Triangles/64K Apple Ii, Iie, Iic/Disk, Backup, Teacher's Guide/Quick Reference Cards/Book No 1340-Fs [online]. —— Cincinnati : Sunburst Communications, 1986. – Available from : http://www.amazon.com/Geometric-Supposer-Triangles-Teachers-Reference/dp/9996313344/ (in English)
11. Greenwald S. B. Computer-assisted explorations in mathematics: pedagogical adaptations across the Atlantic : a report on the development of CATAM at Cambridge University and the Mathematics Project Laboratory at MIT / Suzanne B. Greenwald, Haynes R. Miller // University Collaboration for Innovation: Lessons from the Cambridge-MIT Institute, Sense Publishers. — Massachusetts Institute of Technology, Department of Mathematics, 2007. — P. 121–131. (in English)
12. Harvey J. G. Mathematics Testing with Calculators: Ransoming the Hostages / John G. Harvey // Mathematics assessment and evaluation: Imperatives for mathematics educators / Edited by Tomas A. Romberg. – New York : State University of New York, 1992. – P. 139–168. (in English)
13. Harvey B. A case study: the Lincoln-sudbury regional high school [online] / Brian Harvey // Computer Science Division : EECS at UC Berkeley, 2012. – Available from : http://www.cs.berkeley.edu/~bh/lsrhs.html (in English)
14. Haven G. IME-86: the only electronic desk calculator thet can grow into a compunrt [online] / Grand Haven

// Hope College Department of Mathematics, 1967. – Available from: http://www.math.hope.edu/tanis/History/ime-calculator-purchase.pdf (in English)
15. Horn J. K. The HP-71B «Math machine» 25 years old [online] / Joseph K. Horn, Richard J. Nelson, Dale Thorn // Hewlett-Packard Development Company, [2007]. — Available from: http://h71028.www7.hp.com/enterprise/downloads/The%20HP-71B%20Math%20Machine%20V2e.pdf (in English)
16. John J Wavrik [online] / John J. Wavrik // University of California. —Available from: http://math.ucsd.edu/~jwavrik/(in English)
17. Kernighan B. W. The UNIX programming environment / Brian W. Kernighan, Rob Pike. — Prentice-Hall, 1984. — 357 p.
18. Molnar A. Computers in Education: A Brief History [online] / Andrew Molnar // THE Journal: Technological Horizons in Education. — 06.01.1997. —Available from: http://thejournal.com/Articles/1997/06/01/Computers-in-Education-A-Brief-History.aspx?Page=3(in English)
19. Multipoint data communication system with collision detection : US 4063220 A [online] / David R. Boggs, Butler W. Lampson, Robert M. Metcalfe, Charles P. Thacker. — March 31, 1975. —Available from: http://www.google.com/patents/US4063220(in English)
20. Navarro A. B. Turing Machine Simulator. / Navarro Aaron B. // Computers in mathematics and science teaching, 1981. — v1. — № 2. —P. 25–26. (in English)
21. Pea R. D. Cognitive technologies for mathematics education / Roy D. Pea //Cognitive science and mathematics education, 1987. — C. 89–122. (in English)
22. Sinclair N. Understanding and projecting ICT trends in mathematics education / NathalieSinclair, Nicholas Jackiw// Teaching secondary mathematics with ICT / Edited by Sue Johnston-Wilder and David Pimm. — New York : Open University Press, 2005. — P. 235–253. (in English)
23. Skinner B. F. The technology of teaching / B. F. Skinner. — New York : Appleton-Century-Crofts, 1968. — 271 p. (in English)
24. Suppes P. Computer-based mathematics instruction : The First Year of the Project : (1 September 1963 to 31 August 1964) [online] / Patrick Suppes. // Bulletin of the International Study Group for Mathematics Learnin, 1965 p. — №3. — P.7–22. —Available from: http://suppes-corpus.stanford.edu/article.html?id=54-1(in English)
25. Tanis A. A computer-based laboratory for mathematical statistics and probability [online] / A. Tanis // Computer lab for statistics and probability : Computer Hope. — Available from: http://www.math.hope.edu/tanis/History/ccuc-8%201977.pdf. (in English)
26. Tanis E An Experimental Approach to the Central Limit Theorem [online] / E. Tanis, Deanna Gross // Annual Spring Meeting of the Michigan Section of the Mathematical Association of America : Hope College Department of Mathematics, 1969. —Available from: http://www.math.hope.edu/tanis/History/clt-deanna-1969-rotated.pdf(in English)
27. Tanis E. Hope College : Department of Mathematics [online] / Elliot Tanis // Academic departments math. — Available from: http://www.math.hope.edu/tanis/oldpapersgiven.html (in English)
28. Tanis E. A. Concepts in probability and statistics illustrated with the computer / Elliot A. Tanis // Michigan association of computer users for learning. — 1978. —Vol II. — № 1. — P. 6–17.
29. Tanis E. A. Laboratory manual for probability and statistical inference / Elliot A. Tanis —Iowa : CONDUIT, 1977. — 88 p. (in English)
30. Tanis E. A. Theory of Probability and statistics illustrated by the computer [online] / Elliot A. Tanis // GLCA Computing symposium at Wabash University : Hope College Department of Mathematics. — March 7–8. — 1972. – Available from: http://www.math.hope.edu/tanis/History/wabash%201972-rotated.pdf (in English)
31. The mathematical sciences : Undergraduate education : a report // National research council (U.S.). Panel on undergraduate education in mathematics. — National Academies, 1968. — 113 p. (in English)
32. This month in ed tech history archives [online] // The Software & Information Industry Association, 2013. . — Available from: http://www.siia.net/index.php?option=com_content&view=article&id=314:this-month-archive&catid=159:education-articles&Itemid=328 (in English)
33. UNIX BSD information [online] // Computer Hope. – Available from: http://www.computerhope.com/unix/bsd.htm (in English)
34. Usluel Y. K. A Structural Equation Model for ICT Usage in Higher Education / YaseminKocakUsluel, PetekAskar, Turgay Bas // Educational Technology & Society, 2008. — № 11 (2). — P. 262–273. (in English)
35. Vleck T. V. Dennis M. Ritchie — A. M. Turing Award Winner [online] / Tom Van Vleck // Association for Computing Machinery, 2012. — Available from: http://amturing.acm.org/award_winners/ritchie_1506389.cfm (in English)

36. Watson P. G. Using the computer in education: A briefing for school decision makers / Paul G. Watson. — New Jersey : Educational Technology, 1972. — 128 p. (in English)
37. Wavrik J. J. Computers and the Multiplicity of Polynomial Roots / John J. Wavrik // The American Mathematical Monthly : An official journal of the Mathematical Association of America. — 1982. — Vol. 89, No. 1. — P. 34–36, 45–56. (in English)
38. World University Rankings 2010-2011 — Times Higher Education [online] // TSL Education Ltd. — London, 2012. — Available from:http://www.timeshighereducation.co.uk/world-university-rankings/2010-11/world-ranking. (in English)